\documentclass[12pt]{article}

\usepackage{amssymb,amsmath,amsfonts,graphicx,subfigure,color,multicol}
\usepackage[left=2cm, right=2cm, top=2cm]{geometry}
\usepackage{empheq}
\usepackage{amssymb}
\usepackage{mathtools}
\usepackage{bbm}
\usepackage{dsfont}
\usepackage{epstopdf}
\usepackage{hyperref}
\usepackage{algpseudocode}
\usepackage{algorithm}

\newtheorem{thm}{Theorem}[section]

\newtheorem{dfn}[thm]{Definition}
\newtheorem{rem}[thm]{Remark}
\newtheorem{exa}[thm]{Example}

\newcommand{\Ss}{{\mathcal{S}}}

\usepackage{hyperref}
\usepackage{amssymb,amsmath,amsfonts,graphicx,color,multicol}
\newcommand{\Cc}{{\mathbb{C}}}
\newcommand{\al}{\alpha}
\newcommand{\Ee}{{\mathcal{E}}}

\begin{document}

\title{On the definition of the stability region of multistep methods}

\author{Lajos L\'oczi\thanks{{\texttt{LLoczi@inf.elte.hu}}, Department of Numerical Analysis, Faculty of Informatics, E\"otv\"os Lor\'and University, Budapest, Hungary\newline
\indent \ \,  }}


\maketitle

\begin{abstract}
The usual definition of the stability region of implicit multistep methods often implies that there are some isolated points of stability within the region of instability of the numerical method. 
These isolated stable points may appear when the leading coefficient of the characteristic polynomial of the  method vanishes---they cannot be detected by the well-known root locus method, and their existence renders many results about stability regions problematic.
It is suggested that the definition of the stability region should exclude such isolated points.
\end{abstract}

\section{Introduction}

The aim of this short note is to point out the presence of certain isolated points of stability  
within the region of instability of some common implicit numerical methods.
We argue that these points should not be included in the definition of the stability region.

Stability properties of a broad class of numerical methods (including 
Runge--Kutta methods, linear multistep methods, or multistep multiderivative methods)
for solving initial value problems of the form 
\begin{equation}\label{IVP}
y'(t)=f(t,y(t)), \quad y(t_0)=y_0
\end{equation} 
can be
analyzed by studying the stability region of the method. When an $s$-stage $k$-step 
method ($s\ge 1$, $k\ge 1$ fixed positive integers) with constant
step size $h>0$ is applied to
the linear test equation 
\[
y'=\lambda y \quad\quad(\lambda\in \Cc \text{ fixed},\  y(0)=y_0 \text{ given}),
\]
 the method yields a numerical solution $y_n$ ($n\in\mathbb{N}:=\{0,1,2,\ldots\}$) that
satisfies a recurrence relation of the form \cite{jeltschnevanlinna}
\begin{equation}\label{recurrence}
\left\{
\begin{aligned}
   & \sum_{j=0}^s \sum_{\ell=0}^k a_{j, \ell} \, \mu^j \, y_{n+\ell} =0, \quad \ n\in\mathbb{N}, \\ 
   &   a_{j, \ell}  \in \mathbb{R}, \ \  \sum_{j=0}^s |a_{j,k}|>0, \ \  \mu:=h\lambda.
\end{aligned}
\right.
\end{equation}
The characteristic polynomial associated with the method takes the form
\begin{equation}\label{phidef}
\Phi(\zeta,\mu):=   \sum_{j=0}^s \sum_{\ell=0}^k a_{j, \ell} \, \mu^j \, \zeta^\ell \quad
(\zeta\in\Cc, \,\mu\in\Cc).
\end{equation}
The stability region of the method is defined \cite{hairerwanner2} as
\begin{equation}\label{stabregdef}
\Ss:=\{ \mu\in\mathbb{C} : \text{all roots } \zeta_m(\mu) \text{ of } \zeta \mapsto \Phi(\zeta,\mu)
 \text{ satisfy } |\zeta_m(\mu)|\le 1,
\end{equation}
\[
\text{ and multiple roots satisfy } |\zeta_m(\mu)|<1\}.
\]
\begin{rem}
In \cite{jeltschnevanlinna}, we have $\mu\in\overline\Cc$  in \eqref{stabregdef} instead of 
$\mu\in\Cc$.
\end{rem}
The above definition \eqref{stabregdef} characterizes the boundedness of the sequence $y_n$ ($n\in\mathbb{N}$) generated by
the numerical method \eqref{recurrence} for any \emph{possible} set of initial values $y_0, y_1, \dots, y_{k-1}$ and step size $h>0$.

\begin{exa}\label{firstexample}
A linear $k$-step method \cite{hairerwanner1,hairerwanner2} approximating the solution of the initial value problem \eqref{IVP} can be written as
\begin{equation}\label{multistepdef}
\sum_{\ell=0}^k (\alpha_\ell y_{n+\ell}-h \beta_\ell f_{n+\ell})=0,
\end{equation}
where the numbers $\al_\ell\in\mathbb{R}$ and $\beta_\ell\in \mathbb{R}$ ($\ell=0, \ldots, k$)  are the method coefficients, $\al_k\ne 0$, $t_m$ is defined as $t_0+m h$ ($m\in\mathbb{N}$), and $f_m$ stands for $f(t_m, y_m)$.
The numerical solution $y_n$ approximates the exact solution $y$ at time $t_n$.  
For $k=1$ we have a one-step method, while for $k\ge 2$ the scheme is called a multistep method. 
The method is implicit, if $\beta_k\ne 0$. By setting
\[\varrho(\zeta):=\sum_{\ell=0}^k \al_\ell \zeta^\ell\quad \text{and}\quad \sigma(\zeta):=\sum_{\ell=0}^k \beta_\ell \zeta^\ell,\]
the associated characteristic polynomial \eqref{phidef} is 
$\Phi(\zeta,\mu)=\varrho(\zeta)-\mu \sigma(\zeta).
$\end{exa}
\begin{exa}\label{secondexample}
Multiderivative multistep methods (or generalized multistep methods) 
extend the above class of methods by evaluating the derivatives of $f$ 
at certain points as well.
For example, a second-derivative $k$-step method \cite{hairerwanner2} has the form 
\[
\sum_{\ell=0}^k (\alpha_\ell y_{n+\ell}-h \beta_\ell f_{n+\ell}-h^2 \gamma_\ell g_{n+\ell})=0,
\]
where $g_m:=g(t_m,y_m)$ with $g(t,y):=\partial_1 f(t,y)+\partial_2 f(t,y)\cdot f(t,y)$, and 
the method is determined by the real coefficients $\al_\ell$ ($\al_k\ne 0$), $\beta_\ell$ and $\gamma_\ell$.
The associated characteristic polynomial \eqref{phidef} is now
$
\Phi(\zeta,\mu)=\sum_{\ell=0}^k (\al_\ell-\mu \beta_\ell-\mu^2 \gamma_\ell)\zeta^\ell.
$
\end{exa}

\section{Vanishing leading coefficient of the characteristic polynomial}

The characteristic polynomial of the implicit Euler method with $s=k=1$ is
$\Phi(\zeta,\mu)=(1-\mu)\zeta-1$.
We have $1\in\Ss$, because $\Phi(\zeta,1)=0$ has no roots in $\mathbb{C}$, so \eqref{stabregdef} is satisfied vacuously. For $\mu\ne 1$,
$\Phi(\zeta,\mu)=0$ if and only if $\zeta=1/(1-\mu)$. Hence
\begin{equation}\label{exampleimplicitEuler}
\Ss=\{\mu\in\Cc : |\mu-1|\ge 1\}\cup\Ee
\end{equation}
with $\Ee=\{1\}$. In particular, $1\in\partial \Ss$, the boundary of $\Ss$.

Motivated by the above example, let us rewrite $\Phi$ in \eqref{phidef} as
$
\Phi(\zeta,\mu)=: \sum_{\ell=0}^k C_\ell(\mu)\zeta^\ell
$
with suitable polynomials $C_\ell$. The leading coefficient $C_k$ does not vanish
identically because of the assumption $\sum_{j=0}^s |a_{j,k}|>0$ in
\eqref{recurrence}, or $\al_k\ne 0$ in Examples \ref{firstexample} and \ref{secondexample}.
For implicit methods, $C_k$  is a polynomial of degree
at least 1, so the finite set 
\begin{equation}\label{Edef}
\Ee:=\{\mu\in\Cc : C_k(\mu)=0\}
\end{equation}
is non-empty.

\begin{figure}[H]
\subfigure{
\includegraphics[width=0.45\textwidth]{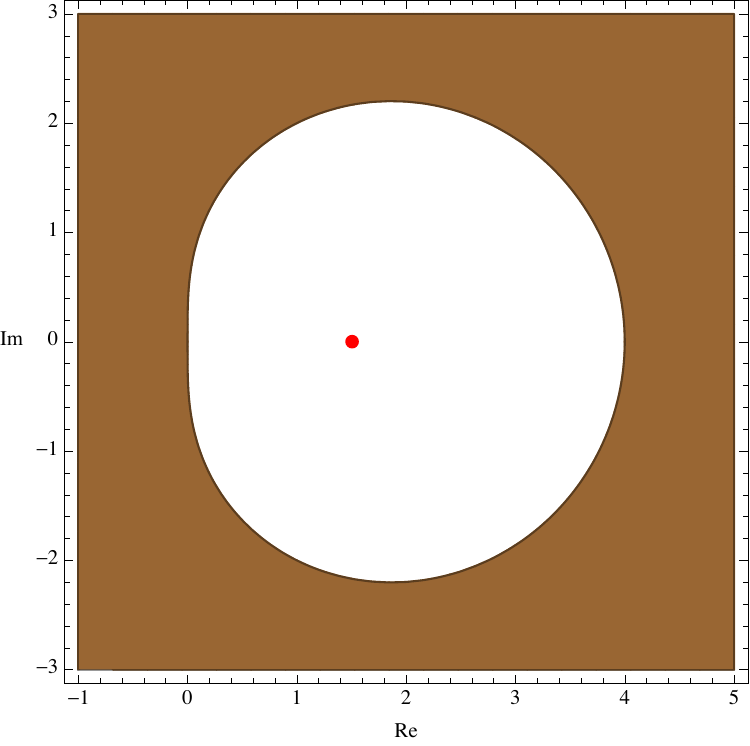}}
\subfigure{
\includegraphics[width=0.45\textwidth]{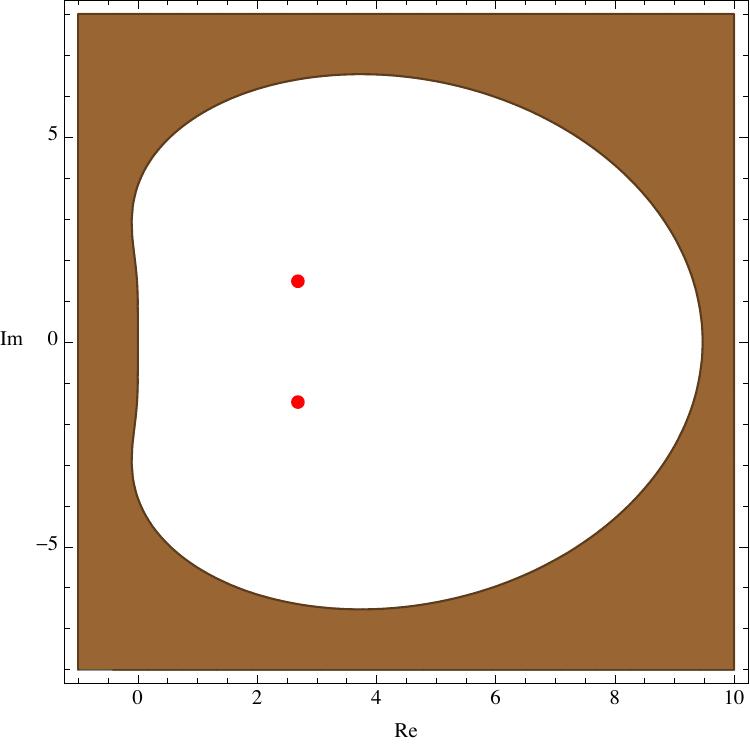}}
\caption{The \emph{left} figure shows the stability region of the $2$-step BDF method in brown and red according to definition \eqref{stabregdef}. The red dot is the unique element of $\Ss\cap\Ee=\{3/2\}$. 
The \emph{right} figure shows the stability region of the $3$-step Enright method 
(member of the family presented in Example \ref{secondexample}) in brown and red. 
For this method we have 
$\Phi(\zeta,\mu)=\left(\frac{19 \mu ^2}{180}-\frac{307 \mu }{540}+1\right)\zeta ^3 +
\left(-\frac{19 \mu }{40}-1\right)\zeta ^2+\frac{\mu }{20}\zeta -\frac{7 \mu }{1080}$, so 
$\Ee=\left\{ \left(307\pm i \sqrt{28871}\right)/114\right\}$ (represented by the two red dots)
and $\Ee\subset\Ss$. 
\label{fig}}
\end{figure}

Besides the implicit Euler method, there are many examples of classical implicit numerical methods
when all the complex roots of the polynomial $\Phi(\cdot,\mu^*)$ have modulus
strictly less than $1$ for some $\mu^*\in \Ee$, hence $\mu^*\in\Ss$.
\begin{exa} The characteristic polynomial of the $2$-step BDF method  \cite{hairerwanner2}  is
$\Phi(\zeta,\mu)=(3-2\mu) \zeta ^2-4 \zeta +1$. Its stability region is depicted in
Figure \ref{fig}. Now $\Ee=\{3/2\}\subset\Ss$, because the 
unique root of $\Phi(\zeta,3/2)=0$ is $\zeta=1/4$.
\end{exa}
\begin{exa}
For several other BDF,
implicit Adams, or 
Enright methods \cite{hairerwanner2} (see Figure \ref{fig}) we have the inclusion $\varnothing\ne\Ee\subset\Ss$.
\end{exa}

Now we point out some consequences of the definition \eqref{stabregdef}.\\


\noindent \textbf{Observation 1.}
If the step size $h>0$ of the method \eqref{recurrence} is chosen in a way that $\mu=h\lambda\in\Ee$,
then the order of the recurrence becomes strictly less than $k$, hence, in general, the initial values $y_0, y_1, \dots, y_{k-1}$ cannot be chosen arbitrarily.\\

\noindent \textbf{Observation 2.} Recursions
with almost zero leading coefficients can be highly unstable with respect to small perturbations.
This renders the corresponding numerical method 
 useless in practice. For example, let us consider the recursion 
corresponding to the $2$-step BDF method
$(3-2\mu)y_{n+2}-4y_{n+1}+y_n=0$. For $\mu=3/2\in\Ss\cap\Ee$, $\displaystyle \lim_{n\to +\infty}y_n=0$ for any starting value $y_0$, but
for small $\varepsilon>0$ and $0<|\mu-3/2|<\varepsilon$, the sequence $|y_n|$ quickly ``blows up"
for generic starting values, since the absolute value of one root of the characteristic polynomial 
$(3-2\mu)\zeta^2-4\zeta+1=0$ is large.\\

\noindent \textbf{Observation 3.}
One way to study $\Ss$---or, more precisely, $\partial \Ss$ (the boundary of $\Ss$)---in the complex plane is to plot the 
\emph{root locus curve} corresponding to the method \cite{hairerwanner2}. 

For methods in Example \ref{firstexample}, $\Phi$ is linear in $\mu$, so $\Phi(\zeta,\mu)=0$ implies $\mu=\varrho(\zeta)/\sigma(\zeta)$ (for $\sigma(\zeta)\ne 0)$. The
root locus curve is then the parametric curve 
\begin{equation}\label{muasaratio}
[0,2\pi)\ni \vartheta\mapsto \mu(\vartheta):=\frac{\varrho\left(e^{i\vartheta}\right)}{\sigma
\left(e^{i\vartheta}\right)}.
\end{equation}
For methods in Example \ref{secondexample}, the
equation $\Phi\left(e^{i\vartheta},\mu\right)=0$ is quadratic in $\mu$ and can be solved to obtain two root locus curves 
\begin{equation}\label{mu12asaquadratic}
[0,2\pi)\ni\vartheta\mapsto \mu_{1,2}(\vartheta)
\end{equation}
corresponding to the method. In general, the root locus curve of the method \eqref{recurrence} is defined 
in \cite[Definition (2.21)]{jeltschnevanlinna} as
\[
\Gamma:=\{ \mu\in \overline{\Cc} : \exists \zeta\text{ with } |\zeta|=1 \text{ and } 
\Phi(\zeta,\mu)=0\}.
\]

Simple examples show that the root locus curve $\Gamma$ can be a \emph{proper subset} of the boundary of the stability region $\partial \Ss$. In  \cite[Corollary 2.6]{jeltschnevanlinna} it is shown however that for methods satisfying 
Property C (see \cite[Definition 4.7]{hairerwanner2} or \cite[Formula (2.9)]{jeltschnevanlinna}), one has
$\partial\Ss =\Gamma$. 

According to \cite[Section V.4]{hairerwanner2}, all one-step methods have Property C, so, for example, the implicit Euler method also has. And indeed, applying  \cite[Proposition 2.7]{jeltschnevanlinna} to 
the implicit Euler method we get that the polynomials $\varrho(\zeta)=\zeta-1$ and $\sigma(\zeta)=\zeta$ have no common roots and 
$\varrho/\sigma$ is univalent on the set $\{ z\in\overline\Cc : |z-1|\ge 1\}$, so $Q(\mu)=1/(1-\mu)$ has Property C, thus
$\partial\Ss =\Gamma$. Since now $\Phi(\zeta,1)=\varrho(\zeta)-\sigma(\zeta)=-1$, we see that 
$1\notin\Gamma=\partial\Ss$. On the other hand, we have seen in \eqref{exampleimplicitEuler} that $1\in \partial\Ss$ due to definition \eqref{stabregdef}. 
This apparent contradiction seems to indicate that the authors of \cite{jeltschnevanlinna} interpreted 
definition \eqref{stabregdef} \emph{intuitively}: a root $\zeta=\infty$ is tacitly introduced
as soon as the leading coefficient $C_k(\mu)$ becomes zero. So \cite[Corollary 2.6]{jeltschnevanlinna},
for example, actually relies on Definition \ref{fulldegree} below, rather than on definition \eqref{stabregdef}.
\begin{rem}
In Figure \ref{fig}, elements of the set $\Ss\cap\Ee$ (the red dots) are the isolated elements of $\partial\Ss$,  
and are \emph{not} part of the corresponding root locus curves. We remark that there are examples where an isolated element of $\partial\Ss$  is found on the root locus curve.
\end{rem}

\section{Conclusion}

Based on the above observations it seems reasonable to refine the definition of the stability region of multistep methods as follows (affecting only the class of implicit methods). 
\begin{dfn}\label{fulldegree}
The stability region of a linear multistep or multiderivative multistep method 
with $k\ge 1$ step(s) and with stability polynomial 
\eqref{phidef} is defined as
\begin{gather}
\Ss:=\{ \mu\in\mathbb{C} : \text{the degree of } \Phi(\cdot,\mu) \text{ is exactly }
k,\nonumber\\
\text{ all roots } \zeta_m(\mu) \text{ of } \zeta \mapsto \Phi(\zeta,\mu)
 \text{ satisfy } |\zeta_m(\mu)|\le 1, 
\text{ and multiple roots satisfy } |\zeta_m(\mu)|<1\}.\nonumber
\end{gather}
\end{dfn}
\medskip
\begin{rem}
Definition \ref{fulldegree} with the non-vanishing leading coefficient
essentially appears, for example, in \cite[Section 2.1]{spijker2013} 
(where it is formulated for linear multistep methods, that is, for $s=1$ in \eqref{recurrence}), or in \cite[Section 2]{spijker2017}. 
\end{rem}
\bibliographystyle{spmpsci}      

\begin{thebibliography}{}

\bibitem{jeltschnevanlinna}
R.~Jeltsch, O.~Nevanlinna, \emph{Stability and Accuracy of Time Discretizations for Initial Value Problems}, 
Numer. Math., Vol.~40, 245--296 (1982)

\bibitem{hairerwanner1}
E.~Hairer, S.~N{\o}rsett, G.~Wanner, \emph{Solving Ordinary Differential Equations I. Nonstiff
Problems}, Springer, Berlin (2009)

\bibitem{hairerwanner2}
E.~Hairer, G.~Wanner, \emph{Solving Ordinary Differential Equations II. Stiff
and Differential-Algebraic Problems}, Springer, Berlin (2002)

\bibitem{spijker2013}
M.~N.~Spijker, \emph{The existence of stepsize-coefficients for boundedness of linear multistep
methods},  Appl. Numer. Math.,  Vol.~63, 45--57 (2013) 

\bibitem{spijker2017}
M.~N.~Spijker, \emph{Stability and boundedness in the numerical solution of initial value problems},  Math. Comp., Vol.~86, No. 308, 2777--2798  (2017)


%
%
\end{thebibliography}

{\footnotesize

}

\end{document}